\documentclass[letter,12pt]{report}

\usepackage{graphicx}    
\usepackage{tikz}
\usepackage{float}
\usepackage{faktor}
\usepackage{subcaption}
\usepackage[centertableaux]{ytableau}
\usepackage{caption}
\usepackage{colortbl}
\usepackage{comment}
\usepackage[utf8]{inputenc}

\usepackage{listings}
\usepackage{xcolor}

\definecolor{codegreen}{rgb}{0,0.6,0}
\definecolor{codegray}{rgb}{0.5,0.5,0.5}
\definecolor{codepurple}{rgb}{0.58,0,0.82}
\definecolor{backcolour}{rgb}{0.95,0.95,0.92}

\lstdefinestyle{mystyle}{
	backgroundcolor=\color{backcolour},   commentstyle=\color{codegreen},
	keywordstyle=\color{magenta},
	numberstyle=\tiny\color{codegray},
	stringstyle=\color{codepurple},
	basicstyle=\ttfamily\footnotesize,
	breakatwhitespace=false,         
	breaklines=true,                 
	captionpos=b,                    
	keepspaces=true,                 
	numbers=left,                    
	numbersep=5pt,                  
	showspaces=false,                
	showstringspaces=false,
	showtabs=false,                  
	tabsize=4
}

\newcommand{\svdots}{%
	\vbox{
		\scriptsize \baselineskip 2pt \lineskiplimit 0pt
		\hbox {.}\hbox {.}\hbox {.}\kern-0.75pt
	}%
}

\usepackage{setspace} 
\usepackage{hyperref}
\usepackage{mathrsfs}
\usepackage{calc}
\usepackage{fancyhdr}
\usepackage{amsmath}
\usepackage{amssymb}
\usepackage[nottoc,notlof,notlot]{tocbibind} 
\usepackage{lipsum} 
\usepackage[top=1in, bottom=1in, left=1in, right=1in]{geometry} 
\usepackage[titles]{tocloft}

\usepackage{amsthm}
\newtheorem{thm}{Theorem}
\newtheorem*{thm*}{Theorem}
\newtheorem{defn}{Definition}

\usepackage{titlesec}

\titleformat{\chapter}{\centering\normalfont\large}{\thechapter}{1em}{} 
\titleformat{\section}{\centering\normalfont\large}{\thesection}{1em}{}
\titleformat{\subsection}{\centering\normalfont\large}{\thesubsection}{1em}{}

		\titlespacing*{\chapter}{0pt}{0.25in}{20pt}
		\fancypagestyle{sxoli2}{%
\fancyhf{}
\pagenumbering{arabic}
\rhead{\thepage}

        \fancypagestyle{plain}{%
        \rhead{\thepage}
        }
}
\usepackage{etoolbox}
\makeatletter
\patchcmd{\chapter}{\if@openright\cleardoublepage\else\clearpage\fi}{}{}{}
\makeatother

\begin{document}
		
    \pagenumbering{arabic}
		\begin{center}
		\singlespacing
		\textbf{ALGEBRAIC RELATIONS VIA A MONTE CARLO SIMULATION}\\
		\doublespacing
		\footnotesize{ALISON ELAINE BECKER}\\
		\singlespacing
		\end{center}
		\begin{flushleft}
		\singlespacing
\footnotesize{ABSTRACT. The conjugation action of the complex orthogonal group on the polynomial functions on $n \times n$ matrices gives rise to a graded algebra of invariant polynomials. A spanning set of this algebra is in bijective correspondence to a set of unlabeled, cyclic graphs with directed edges equivalent under dihedral symmetries. When the degree of the invariants is $n+1$, we show that the dimension of the space of relations between the invariants grows linearly in $n$. Furthermore, we present two methods to obtain a basis of the space of relations. First, we construct a basis using an idempotent of the group algebra referred to as Young symmetrizers, but this quickly becomes computationally expensive as $n$ increases. Thus, we propose a more computationally efficient method for this problem by repeatedly generating random matrices using a Monte Carlo algorithm.}
		\end{flushleft}
	\singlespacing
    \chapter{Introduction}
The goal of this work is to understand the relations between invariants of the conjugation action of the complex orthogonal group, $O_n(\mathbb{C})$, on the polynomial functions on $n \times n$ matrices, $\mathcal{P}(M_n)$. We present results about the dimension of the space of relations, denoted $\mathcal{R}\mathcal{E}\mathcal{L}$, and further discuss two methods for obtaining a basis of this space. First, we construct a basis of relations by using Young symmetrizers, however, this can be computationally expensive. Thus we propose a more efficient method to determine a basis of $\mathcal{R}\mathcal{E}\mathcal{L}$ by repeatedly generating random matrices via a Monte Carlo algorithm. 

Throughout this paper we let $GL_n$ denote the complex general linear group, and we let $O_n$ denote the complex orthogonal group. 
In Section $3$ we discuss our main results about the dimension of $\mathcal{R}\mathcal{E}\mathcal{L}$ when the degree of the polynomial invariants is equal to $n+1$. Again, these are the invariants under the conjugation action of $O_n$ on $\mathcal{P}(M_n)$,	
	\begin{equation*}
		g\cdot f(x) = f(g^{T}xg)
	\end{equation*}
where $x \in M_n, f \in \mathcal{P}(M_n)$ and $g \in O_n$. 

We have that $\mathcal{P}^d(M_n)$, the homogeneous polynomials of degree $d$, are a finite dimensional representation of $O_n$. Therefore, we have a graded structure on $\mathcal{P}(M_n)$, 
	\begin{equation*}
		\mathcal{P}(M_n) = \bigoplus\limits_{d \geq 0}\mathcal{P}^d(M_n)
	\end{equation*}
and thus a graded algebra of invariant polynomials, $\mathcal{P}(M_n)^{O_n}$. It is shown in \cite{Procesi} that this algebra is generated by functions of the form $Tr(x^{a_1}(x^T)^{a_2}x^{a_3}(x^T)^{a_4}\cdots x^{a_{m-1}}(x^T)^{a_m})$ for $x\in M_n$. Products of these functions span the set of invariants, however, they are not linearly independent if the degree of the polynomials is greater than $n$. 

There are no relations between invariants in what we define as the \textit{stable range}, where the degree of the invariant polynomials is smaller than $n$ \cite{Jeb}. The first time relations arise is just outside of the stable range, where the degree of the polynomials is $n+1$. The following theorem gives the dimension of the space of relations between the polynomial invariants, $\mathcal{R}\mathcal{E}\mathcal{L}_{n+1}$, which we see grows linearly in $n$.

\begin{thm*}
	Let $n$ be a positive integer. The dimension of the space of relations, $\mathcal{R}\mathcal{E}\mathcal{L}_{n+1}$, between the degree $n+1$ invariants of the $O_n$ conjugation action on $\mathcal{P}(M_n)$ is equal to  
	\begin{equation*}
	dim(\mathcal{R}\mathcal{E}\mathcal{L}_{n+1}) = \begin{cases}
	\dfrac{n}{2} +1 &n \textnormal{  even} \vspace{.3cm}\\ 
	\dfrac{n+3}{2} &n \textnormal{  odd} \end{cases}
	\end{equation*}
\end{thm*} 

In order to prove Theorem \ref{my result}, we rely on correspondences between the elements of $\mathcal{P}(M_n)^{O_n}$ and several different spaces. In Section $2$ we discuss a useful bijection from \cite{Jeb} between the invariant polynomials and Necklace diagrams, which are unlabeled, cyclic graphs embedded in the plane. 

Next, we present two ways to determine a basis of $\mathcal{R}\mathcal{E}\mathcal{L}$. In Section \ref{Constructing a Basis} we explicitly describe the relations using idempotents of the group algebra, $\mathbb{C}[S_n]$, called Young symmetrizers. However, the size of these elements grows exponentially in $n$, and thus it becomes computationally overwhelming to obtain a basis using this method as $n$ increases. 

Therefore, in Chapter \ref{Monte Carlo}, we propose a new approach to explicitly define a basis of the relations that is more efficient and does not rely on Young symmetrizers. This algorithm is designed using a Monte Carlo simulation that repeatedly generates random matrices to give numerical values to the invariants, and we solve a linear system to recover a basis of relations. 

All of our code is written in Python and Sage \cite{Sage} and run on the SageMath cloud servers.

\chapter{Invariants and correspondences}

\textbf{$\mathbf{2.1}$ $\mathbf{O_n}$ invariant polynomials on matrices.}\label{invariants} Consider the complex general linear group, $GL_n$, and the polynomial functions on $n \times n $ matrices, $\mathcal{P}(M_n)$. There is a conjugation action of $GL_n$ on $\mathcal{P}(M_n) $ defined by 
\begin{center} $g\cdot f(x) = f(g^{-1}xg)$ \end{center} 
for $g \in GL_n, x \in M_n,$ and $f \in \mathcal{P}(M_n)$. 

We can see that $Tr(x)$ is invariant under this action, since 
\begin{equation*}
	g \cdot Tr(x) = Tr(g^{-1}xg) = Tr(gg^{-1}x) = Tr(x) 
\end{equation*}

It can be shown \cite{Procesi} that $Tr(x^k)$ is invariant for $k \in \mathbb{N}, k\leq n$, and these functions generate the invariant algebra $\mathcal{P}(M_n)^{GL_n}$.  Furthermore, there are no relations between these invariants. Thus, we restrict the adjoint action of $GL_n$ to the action of a subgroup of $GL_n$, the complex orthogonal group $O_n$, on $\mathcal{P}(M_n)$ defined by 
	\begin{equation*}\label{On action}
		g\cdot f(x) = f(g^{T}xg)
	\end{equation*}  
The homogeneous polynomials of degree $d$, denoted $\mathcal{P}^d(M_n)$, are a finite dimensional representation of $O_n$. Therefore, we have a graded algebra structure, 
\begin{equation*}
\mathcal{P}(M_n) = \bigoplus\limits_{d} \mathcal{P}^d(M_n)
\end{equation*}
and a similar grading of the algebra of polynomial invariants, $\mathcal{P}(M_n)^{O_n}$. 

The polynomial invariants under the $GL_n$ action are also invariant under the orthogonal group, but there are additional invariants under the conjugation action of $O_n$. It is shown in \cite{Procesi}, that $\mathcal{P}(M_n)^{O_n}$ is generated by traces of monomials in $x$ and $x^T$, 
\begin{equation*}\label{traces} 
Tr(x^{a_1}(x^T)^{a_2}x^{a_3}(x^T)^{a_4}\cdots x^{a_{m-1}}(x^T)^{a_m}) 
\end{equation*}
for $x \in M_n, a_i \in \mathbb{Z}^+$. 

Products of these polynomial functions span the set of invariants. When the degree of the invariants is less than or equal to $n$, there are no relations between the polynomials. This is called the stable range, defined below. 

\begin{defn}\label{stable range}
	For $d,n \in \mathbb{N}$, we define the ordered pair $(d,n)$ to be in the \textbf{stable range} if $d \leq n$.
\end{defn}
If the degree, $d$, of the polynomials is greater than the dimension of the defining representation, the spanning set of invariants is not linearly independent. We look at the first occurrence of relations between invariants, which happens just outside of the stable range where the invariant polynomials have degree $d=n+1$. There is a useful correspondence between these polynomials and unlabeled cyclic graphs called Necklace diagrams.
\vspace{.5cm}
\newline 
\textbf{$\mathbf{2.2}$ Necklace diagrams.}\label{Necklace diagrams} There is a bijective correspondence between unlabeled cyclic graphs with oriented edges and the polynomial generators of the invariant ring $\mathcal{P}(M_n)^{O_n}$ detailed in \cite{Jeb}. These graphs are called Necklace diagrams, defined as follows. 
	\begin{defn}\label{Necklace Diagram defn}
		Let $m \in \mathbb{N}$. We embed an unlabeled, directed, cyclic graph with $m$ nodes in the plane and centered at the origin; each edge is given an orientation of clockwise or counterclockwise determined by choosing an arbitrary edge, $E$, and noting the direction traveled to subsequent edges. 
		
		These graphs, called \textbf{Necklace diagrams} and denoted $N_m$, have the following structure: 
			\begin{itemize}
				\item An edge may join a node to itself
				\item At most two edges may join two different nodes	
			\end{itemize}
		Furthermore, diagrams are considered equivalent under dihedral symmetries of rotation and reflection.	
	\end{defn}
Let $\mathscr{D}_d$ denote the set of all (not necessarily connected) Necklace diagrams in which each connected component is a Necklace diagram $N_{m_i}$, and there are $d$ total nodes in the diagram. That is, define 
	\begin{equation*}
		\mathscr{D}_d = \{N_{m_i} | \hspace{.2cm} N_{m_i} \textnormal{ is a Necklace diagram with $m_i$ nodes, and} \sum_{i}m_i = d \}.
	\end{equation*}
\textbf{Remark.} The combinatorics presented here regarding the construction of Necklace diagrams is also seen in \cite{Brucks}. Here, we allow a diagram with $d$ nodes to consist of possibly disconnected components. In \cite{Brucks}, the author considers the number of connected diagrams with no proper subpattern.   

In the stable range where $d \leq n$, the dimension of the space of invariants under the $O_n$ conjugation action is exactly equal to the number of Necklace diagrams with $d$ total nodes, that is, 
$$|\mathscr{D}_d| = \dim\mathcal{P}^d(M_n)^{O_n}.$$ As a result, there are no relations between the invariants in this setting. Thus, we present a new result about the dimension of the invariant space outside of the stable range, where relations between the invariants arise. To aid in this discussion, it is useful to note a correspondence between Necklace diagrams and a set of fixed-point free involutions. 

By definition, the set of all fixed-point free involutions on a set $A$ is the set of all transpositions of elements in $A$. Let $I_n$ denote the set of all fixed-point free involutions on the set $\{1, 2, \dots, n\}$. Then by \cite{Jeb}, we have a bijective correspondence 
\begin{equation*}\label{Necklace FFI correspondence}
\Theta: \mathscr{D}_d \longrightarrow I_{2d}
\end{equation*}
between the set of Necklace diagrams with $d$ total nodes and involutions without fixed points on $\{1, 2, \dots, 2d \}.$

Thus, we are free to consider the degree $d$ elements of $\mathcal{P}(M_n)^{O_n}$ as products of traces of $x$ and $x^T$, or as Necklace diagrams, or as fixed-point free involutions of the set $\{1,2,\dots,2d \}$. In the next section, we describe an important correspondence between these fixed-point free involutions, $I_{2d}$, and a set of double cosets of the symmetric group.  
\vspace{.5cm}
\newline
\textbf{2.3 Fixed-point free involutions and the symmetric group.} Consider the symmetric group $S_{2n}$, and note the following inclusion: 
	\begin{equation}
		\Delta S_{n} \subseteq S_n \times S_n \hookrightarrow S_{2n}
	\end{equation}
where $\Delta S_{n}$ denotes the diagonally embedded copy of $S_n$ in $S_n \times S_n$, that is, 
	\begin{equation*}\label{diagonal embed}
		\Delta S_{n}= \{(\sigma, \sigma) | \hspace{.1cm} \sigma \in S_n \}
	\end{equation*}
Additionally, we consider $S_n \times S_n$ where the first copy of $S_n$ is composed of permutations of the set $\{1,2, \dots, n \}$ and the second copy of $S_n$ is permutations of the set $\{n+1, n+2, \dots, 2n \}$. Then we have an embedding of $S_n \times S_n$ into $S_{2n}$, and thus we can consider $\Delta S_n$ and $S_n \times S_n$ as subgroups of $S_{2n}$.
	
We define the following element of $S_{2n}$ using disjoint cycle notation: 
$$ \tau = (1 \hspace{.2cm} 2)(3 \hspace{.2cm} 4) \cdots (i \hspace{.2cm} i+1) \cdots (2n-1 \hspace{.2cm} 2n)	$$
and we define $H_n$ as the centralizer of $\tau$ in the group $S_{2n}$, 
$$ H_n = \{\sigma \in S_{2n} | \hspace{.1cm} \sigma\tau = \tau\sigma  \}.  $$
Clearly, by definition $\tau$ is also an element of $I_{2n}$. Now, we can construct a set of double cosets of the group $S_{2n}$ as follows: 
	\begin{equation*}\label{double cosets}
		(\Delta S_{n})\backslash S_{2n}/H_n = \{(S_n)\sigma H_n | \hspace{.1cm} \sigma \in S_{2n}  \} 
	\end{equation*}
There is a bijection, \cite{Jeb}, between these double cosets and elements of the invariant algebra $\mathcal{P}(M_n)^{O_n}$. We note that this construction of invariants is used in our code to find the elements of  $\mathcal{P}(M_n)^{O_n}$.

In the following sections our goal is to understand the dimension of the space of relations between the elements of $\mathcal{P}(M_n)^{O_n}$; in order to do this we rely on the combinatorial object called Littlewood-Richardson numbers. 
\vspace{.5cm}
\newline
\textbf{2.4 Littlewood-Richardson Numbers.} There are several equivalent definitions of Littlewood-Richardson numbers, for this work we consider them as the coefficients that arise from inducing a representation from a subgroup of a symmetric group to the entire symmetric group.

\begin{defn}\label{littlewood}
	Let $\mu$ and $\nu$ be partitions of positive integers $n$ and $m$, respectively. Then the module $Y_{\mu} \otimes Y_{\nu}$ is naturally an $(S_{n} \times S_{m})$-module. Thus we can induce to an $S_{n+m}$ representation, 
		\begin{equation}
			Ind_{S_{n}\times S_{m}}^{S_{n+m}} Y_{\mu} \otimes Y_{\nu} = \bigoplus\limits_{\lambda \vdash n+m} c_{\mu\nu}^{\lambda}Y_{\lambda}
		\end{equation} 
	where the coefficients, $c_{\mu\nu}^{\lambda}$ are \textbf{Littlewood-Richardson numbers}. 
\end{defn}

Furthermore, the Littlewood-Richardson rule states that the coefficients $c_{\mu\nu}^{\lambda}$, count the number of skew semi-standard Young tableaux \cite{YoungTab} of shape $\lambda/\mu$ with weight $\nu$, with the additional restriction that the concatenation of the reversed rows is a lattice word. 
 
It can be shown using character theory that $c_{\mu\nu}^{\lambda} \neq 0$ when $\mu, \nu \subseteq \lambda$. That is, the Young diagrams of $\mu$ and $\nu$ must fit inside the Young diagram of $\lambda$ in order to have nonzero coefficients. 

We are now ready to state and prove a theorem about the dimension of the space of relations between the polynomial invariants under the conjugation action of the complex orthogonal group on $\mathcal{P}(M_n)$. 

\chapter{The dimension of the space of relations}
Under the conjugation action of $O_n$ on $\mathcal{P}(M_n)$, we know that there are no relations between the invariants in the stable range, where the degree of the polynomials is less than or equal to $n$. Therefore, we analyze a space outside of the stable range where there are relations. In this section, we present a result about the dimension of the space relations between the invariant polynomials of degree $d =n+1$. We start by giving an example. 

\textbf{Example:} Let $n=2$. We consider invariant polynomials under the action of $O_2$; we know that there are no relations between the invariants if $d \leq 2$. Thus we take one step outside of the stable range and let $d = 3$, where we consider cubic invariant polynomials under the action of $O_2$ on $\mathcal{P}(M_2)$. Here, 
\begin{equation*}
dim(\mathcal{P}^{3}(M_2))^{O_2} = 5
\end{equation*}
The space consists of the following cubic invariants, 
\begin{equation*}
Tr(x^3), \hspace{.2cm} Tr(x^2x^T), \hspace{.2cm} Tr(x^2)Tr(x), \hspace{.2cm} Tr(xx^T)Tr(x), \hspace{.2cm} Tr(x)^3
\end{equation*}
In the language of Section \ref{Necklace diagrams}, these polynomials correspond to the set of Necklace diagrams with $3$ total nodes and oriented edges. 

\begin{figure}[H]
	\begin{center}
		\resizebox{16.5cm}{!}{
			\begin{tikzpicture}
			\tikzstyle{every node} = [circle,fill,inner sep=1.5pt]
			
			\node (1) at (0,0) {};
			\node (2) at ++(-.5,-1) {};
			\node (3) at ++(.5,-1) {};			
			\foreach \from/\to in {2/1, 1/3, 3/2}
			\draw [ ->] (\from) -- (\to);
			
			\node (4) at (1.5,0) {};
			\node (5) at ++(1,-1) {};
			\node (6) at ++(2,-1) {};			
			\foreach \from/\to in {4/5, 5/6, 4/6}
			\draw [ ->] (\from) -- (\to);
			
			\node (7) at (2.5,-.5) {};
			\node (8) at ++(3,-.5) {};
			\node (9) at ++(3.5,-.5) {};
			\draw[->] (7) to [out=70,  in=110, looseness=2] (8);
			\draw[->] (8) to [out=250, in=290, looseness=2] (7);
			\draw[->] (9) to [out=-10, in=140, looseness=15] (9);
			
			\node (10) at (4.5,-.5) {};
			\node (11) at ++(5,-.5) {};
			\node (12) at ++(5.5,-.5) {};
			\draw[->] (10) to [out=70,  in=110, looseness=2] (11);
			\draw[->] (10) to [out=290, in=250, looseness=2] (11);
			\draw[->] (12) to [out=-10, in=140, looseness=15] (12);
			
			\node (13) at (6.5,-.5) {};
			\node (14) at ++(7,-.5) {};
			\node (15) at ++(7.5,-.5) {};
			\draw[->] (13) to [out=-10, in=140, looseness=15] (13);
			\draw[->] (14) to [out=-10, in=140, looseness=15] (14);
			\draw[->] (15) to [out=-10, in=140, looseness=15] (15);
			
			\node[fill=none] (16) at (0,-2.3) {\tiny $Tr(x^3)$}; 
			\node[fill=none] (17) at (1.5,-2.3) {\tiny $Tr(x^2x^T)$};
			\node[fill=none] (18) at (3.3,-2.3) {\tiny $Tr(x^2)Tr(x)$};
			\node[fill=none] (19) at (5.3,-2.3) {\tiny $Tr(xx^T)Tr(x)$};
			\node[fill=none] (20) at (7,-2.3) {\tiny $Tr(x)^3$};
			
			\node[fill=none] (21) at (0,-1.2) {};
			\node[fill=none] (22) at (1.5,-1.2) {};
			\node[fill=none] (23) at (1.5,-1.9) {};
			\node[fill=none] (24) at (5.3,-1.2) {};
			\node[fill=none] (25) at (7,-1.2) {};
			\node[fill=none] (26) at (5.3,-1.9) {};
			\node[fill=none] (27) at (3.3,-1.9) {};
			\node[fill=none] (28) at (3.3,-1.2) {};
			
			\draw[<->] (21) to (16);
			\draw[<->] (22) to (23);
			\draw[<->] (27) to (28);
			\draw[<->] (26) to (24);
			\draw[<->] (25) to (20);			
			\end{tikzpicture}}
	\end{center}
\vspace{-1cm}
	\caption{\footnotesize{Bijection between Necklace diagrams $\mathscr{D}_3$ and elements of $P^3(M_2)^{O_2}$.}}
\end{figure}
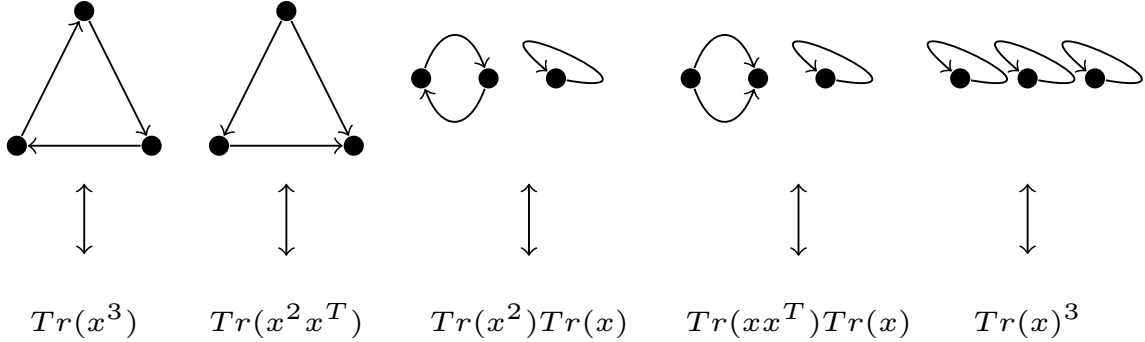
These polynomials span the set of invariants, however, they are not linearly independent. There are actually two relations, 
	\begin{alignat*}{2}
		Tr(x)^3 - 3Tr(x)^2Tr(x) + 2Tr(x^3) &= 0 \\
		2Tr(x^2x^T) - 2Tr(xx^T)Tr(x) + Tr(x)^3 - 3Tr(x^2)Tr(x) &= 0
	\end{alignat*}

The primary result of this section describes how many relations there are between the degree $n+1$ invariants under the action of $O_n$. We start by giving some set up for the main theorem. 
\vspace{.5cm} 
\newline
\textbf{3.1 Dimension of the space of relations between $\mathbf{O_n}$ invariants.} Consider the following inclusion map from the square $n$-dimensional matrices to the square $(n+1)$-dimensional matrices: 
	\begin{alignat*}{2}
		M_n &\hookrightarrow M_{n+1} \\
		x & \hookrightarrow 
		\left( \begin{array}{c|c} 
		x & 0 \\
		\hline
		0 & 0				
		\end{array}  \right)_{n+1}
	\end{alignat*}
Then we have a surjection between the polynomial spaces: 
	\begin{equation*}
		\mathcal{P}[M_{n+1}] \twoheadrightarrow \mathcal{P}[M_{n}]
	\end{equation*}
where we restrict the $(n+1) \times (n+1)$ dimensional matrix down to an $n \times n$ dimensional matrix. Recall, in Section \ref{invariants} we discuss the invariants of this space under the conjugation action of the complex orthogonal group $O_n$. Under this action we have a surjection between the invariant rings: 
	\begin{equation}
		\mathcal{P}[M_{n+1}]^{O_{n+1}} \twoheadrightarrow \mathcal{P}[M_n]^{O_n}
	\end{equation}

\textbf{Remark:} We note that $\mathcal{P}(M_n)$ is a $\mathbb{C}$-algebra of polynomial functions on $M_n$ with a graded structure: 	
	\begin{equation*}
		\mathcal{P}(M_n)= \bigoplus\limits_{d} \mathcal{P}^d(M_n)
	\end{equation*}
where $\mathcal{P}^d(M_n)$ denotes the subspace of homogeneous degree $d$ polynomials, which are a finite dimensional representation of $O_n$. Thus we also have a graded structure on the invariant algebras, 
	\begin{equation*}
		\mathcal{P}(M_n)^{O_n}= \bigoplus\limits_{d} \mathcal{P}^d(M_n)^{O_n}
	\end{equation*}
Therefore if we fix the degree, $d = n+1$, of the invariant polynomials, we have a surjection between the algebras 
\begin{equation}\label{algmap}
\mathcal{P}^{n+1}(M_{n+1})^{O_{n+1}} \twoheadrightarrow \mathcal{P}^{n+1}(M_n)^{{O_n}}
\end{equation}

Note that when $d=1$, we have that $\mathcal{P}^1(M_n) = M_n^*$, which leads to the following identification: 
	\begin{equation*}
		\mathcal{P}(M_n) = \bigoplus\limits_{d}\mathcal{P}^d(M_n) = \bigoplus\limits_{d}\mathcal{S}^d(M_n^*) = \mathcal{S}(M_n^*)
	\end{equation*}
where $\mathcal{S}(M_n^*)$ is the symmetric algebra of the dual space, $M_n^*$, of polynomials on $M_n$.
Then if we again fix the degree $d = n+1$, we see that: 
	\begin{equation*}
		\mathcal{P}^{n+1}(M_n) = \mathcal{S}^{n+1}(M_n^*) = \mathcal{S}^{n+1}(M_n)^*	
	\end{equation*}
Furthermore, as we are working over $\mathbb{C}$, we have that $\mathcal{S}^{n+1}(M_n) \simeq [\otimes^{n+1}(M_n)]^{\Delta S_{n+1}}$. Note that the symmetric tensors are invariant under the natural permutation action of the symmetric group on the tensor factors. Thus map \ref{algmap} can be written as the surjective map: 
	\begin{equation}\label{invariant space map}
		[\otimes^{n+1}(M_{n+1})]^{O_{n+1} \times \Delta S_{n+1}} \twoheadrightarrow [\otimes^{n+1}(M_n)]^{{O_n} \times \Delta S_{n+1}}
	\end{equation}
where $\Delta S_{n+1}$ denotes the diagonally embedded copy of $S_{n+1}$ in $S_{n+1} \times S_{n+1}$. 

Next, we use the following decomposition of $M_n$ in the above map \ref{invariant space map}.  Consider the $Gl_n \times Gl_n$ action on $M_n$ given by:
\begin{equation*}
(g,h)\cdot x = gxh^T
\end{equation*}
for $x\in M_n$ and $(g,h) \in Gl_n \times Gl_n$. Restricting to the diagonal $Gl_n$ action on this space gives: 
\begin{equation*}
(g,g)\cdot x = gxg^T
\end{equation*}
and thus under this action, we have a decomposition of $M_n$ as follows, 
\begin{equation}\label{Mn decomp}
M_n \simeq \mathbb{C}^n \otimes \mathbb{C}^n
\end{equation}
We focus on the $O_n$ decomposition of $M_n$ under the adjoint action of $Gl_n$, where $x \rightarrow gxg^{-1}$. Under this action we have that 
\begin{equation*}
M_n \simeq (\mathbb{C}^n)^* \otimes \mathbb{C}
\end{equation*}
However, consider the map:
\begin{alignat*}{3}
\mathbb{C}^n &\rightarrow &(\mathbb{C}^n)^* \\
v &\mapsto \hspace{.3cm} &\varphi: \mathbb{C}^n &\rightarrow \mathbb{C} \\
&                  &w  &\rightarrow v \cdot w		
\end{alignat*}
where $v, w \in \mathbb{C}^n$ and $v \cdot w$ is the usual dot product. Since the dot product is invariant under the $O_n$ action, we have that as an $O_n$- representation, $\mathbb{C}^n \simeq (\mathbb{C}^n)^*$. Thus, we are free to decompose $M_n$ as in \ref{Mn decomp} using this property that $O_n$ is self-dual. 

\textbf{Remark on the Brauer algebra:} Before we symmetrize and consider the $\Delta S_{n+1}$ action, we can decompose the even dimensional tensor space using the dual, $(\mathbb{C}^n)^*$, which gives the following,
	\begin{equation*}
		[(\otimes^{n+1}\mathbb{C}^n)^* \otimes (\otimes^{n+1}\mathbb{C}^n)]^{O_n} \cong End_{O_n}(\otimes^{n+1}\mathbb{C}^n)
	\end{equation*}
where the endomorphism group is defined to be the Brauer algebra \cite{Goodman}. If we consider the $O_n$ action on the tensor space $\otimes^{k}\mathbb{C}^n$ instead of the traditional general linear group action, then the Brauer algebra replaces the symmetric group in the decomposition of the space via Schur-Weyl duality. A recent discussion can be found in \cite{New}. 

Now, returning to our main thread, using the decomposition of $M_n$ described above, we can write map \ref{invariant space map} as:
	\begin{equation}\label{main map}
		\zeta: [\otimes^{2(n+1)}\mathbb{C}^{n+1}]^{O_{n+1}\times \Delta S_{n+1}} \twoheadrightarrow [\otimes^{2(n+1)}\mathbb{C}^n]^{O_{n} \times \Delta S_{n+1}}
	\end{equation}
We see that in the domain of this map, $[\otimes^{2(n+1)}\mathbb{C}^{n+1}]^{O_{n+1}\times \Delta S_{n+1}}$, we are in the stable range where the degree of the invariants is equal to the dimension of the defining representation. In the codomain, $[\otimes^{2(n+1)}\mathbb{C}^n]^{O_{n} \times \Delta S_{n+1}}$, is where relations arise. Thus, the kernel, $\mathcal{R}\mathcal{E}\mathcal{L}$, of this map 
	\begin{equation*}
		(0) \rightarrow \mathcal{R}\mathcal{E}\mathcal{L} \rightarrow [\otimes^{2(n+1)}\mathbb{C}^{n+1}]^{O_{n+1}\times \Delta S_{n+1}} \rightarrow [\otimes^{2(n+1)}\mathbb{C}^n]^{O_{n} \times \Delta S_{n+1}}
	\end{equation*}
is exactly the space of relations between the degree $n+1$ invariants.

	\begin{defn}\label{rel}
		Let $n \in \mathbb{N}$, and let $\Delta S_{n+1}$ denote the diagonally embedded copy of the symmetric group $S_{n+1}$ in $S_{n+1} \times S_{n+1}$. Then the kernel of the map, 
			\begin{equation*}
				\zeta: [\otimes^{2(n+1)}\mathbb{C}^{n+1}]^{O_{n+1}\times \Delta S_{n+1}} \twoheadrightarrow [\otimes^{2(n+1)}\mathbb{C}^n]^{O_{n} \times \Delta S_{n+1}}
			\end{equation*}
		denoted $\mathcal{R}\mathcal{E}\mathcal{L}$, is defined to be the space of relations between the invariants of the conjugation action of the complex orthogonal group on $\mathcal{P}(M_n)$. 
	\end{defn} 

Our goal is to understand $\mathcal{R}\mathcal{E}\mathcal{L}$ as a vector space in order to compute its dimension, and determine a basis. The following theorem describes the dimension of this vector space when we are just outside of the stable range, that is, when the degree of the invariants is $n+1$. 
	\begin{thm}\label{my result}
	Let $n$ be a positive integer. The dimension of the space of relations, $\mathcal{R}\mathcal{E}\mathcal{L}_{n+1}$, between the degree $n+1$ invariants of the $O_n$ conjugation action on $\mathcal{P}(M_n)$ is equal to  
	\begin{equation*}
		dim(\mathcal{R}\mathcal{E}\mathcal{L}_{n+1}) = \begin{cases}
		\dfrac{n}{2} +1 &n \textnormal{  even} \\
		\dfrac{n+3}{2} &n \textnormal{  odd} \end{cases}
	\end{equation*}
	\end{thm}

In order to determine the dimension of $\mathcal{R}\mathcal{E}\mathcal{L}_{n+1}$, the proof will proceed as follows. First, we consider the projection: 
\begin{equation}
\left[ \otimes^{n+1}(M_n) \right]^{O_n} \longrightarrow \mathcal{P}^{n+1}(M_n)^{O_n}
\end{equation}
where when we project to polynomial space, the tensors become symmetric. The relations in the polynomial algebra pull back into the tensor algebra, and they form an irreducible representation of $S_{2(n+1)}$ corresponding to the partition $[n+1, n+1]$. When this partition is restricted to $S_{n+1} \times S_{n+1}$, there is a multiplicity-free decomposition into irreducible representations $Y^\alpha \otimes Y^ \alpha$ corresponding to size $n+1$ diagrams with at most two parts. Thus, to each $\alpha$ there exists a polynomial relation which corresponds to the $S_{n+1}$ invariant where we embed $S_{n+1}$ into $S_{n+1} \times S_{n+1}$ diagonally. This $\Delta S_{n+1}$ is the symmetric group that symmetrizes to go from the tensor algebra to polynomial space. 

The linear growth of $\mathcal{R}\mathcal{E}\mathcal{L}_{n+1}$ is shown in the highlighted diagonal of the following table: 
\begin{table}[H]
	\centering
	Dimension of $\mathcal{R}\mathcal{E}\mathcal{L}$
	\begin{tabular}{|p{.5cm} || p{1.3cm} p{1.3cm} p{1.3cm} p{1.3cm} p{1.3cm} p{1.3cm} p{1.3cm} p{1.3cm} | }
	\hline
	$d\backslash n$& 1  & 2 & 3 & 4 & 5 & 6 & 7 & 8\\ 
	\hline
	\hline
	1 & 0  & 0  & 0  & 0  & 0  & 0  & 0 & 0 \\
	2 & \cellcolor{blue!20}2  & 0  & 0  & 0  & 0  & 0  & 0 & 0 \\
	3 & 2  & \cellcolor{blue!20} 2  & 0  & 0  & 0  & 0  & 0 & 0 \\
	4 & 5  & 3  & \cellcolor{blue!20}3  & 0  & 0  & 0  & 0 & 0 \\
	5 & 5  & 7  & 4  & \cellcolor{blue!20}3  & 0  & 0  & 0 & 0 \\
	6 & 9  & 13 & 12 & 5  & \cellcolor{blue!20}4  & 0  & 0 & 0 \\
	7 & 9  & 21 & 21 & 14 & 6  & \cellcolor{blue!20} 4  & 0 & 0 \\
	8 & 14 & 33 & 48 & 30 & 19 & 7  & \cellcolor{blue!20} 5 & 0 \\
	9 & 14 & 51 & 75 & 67 & 39 & 21 & 8 & \cellcolor{blue!20} 5 \\
	\hline	
	\end{tabular}

\end{table}

Here, the columns are indexed by $n$, and the rows by the degree of the polynomial invariants in $\mathcal{P}^{n+1}(M_n)^{O_n}$. The data in this table again shows that there are no relations when we are in the stable range where the degree of the invariants is less than or equal to $n$, and it illustrates the linear behavior of the dimension of the space of relations. The values in this table are all generated using code we have written to compute the number of invariants in each scenario.

\textbf{Proof of Theorem \ref{my result}} : Recall we are determining the dimension of the kernel, $\mathcal{R}\mathcal{E}\mathcal{L}_{n+1}$, 
	\begin{equation*}
		\mathcal{R}\mathcal{E}\mathcal{L}_{n+1} \rightarrow [\otimes^{2(n+1)}\mathbb{C}^{n+1}]^{O_{n+1} \times \Delta S_{n+1}}
	\end{equation*}
We have the following decomposition into irreducible representations via Schur-Weyl Duality: 
	\begin{equation}
		\otimes^{2(n+1)}\mathbb{C}^{n+1} \simeq \bigoplus\limits_{\substack{\mu: \hspace{.05cm}\mu \vdash 2(n+1) \\{\ell(\mu)\leq n+1} }} F_{n+1}^{\mu} \otimes Y_{2(n+1)}^{\mu}
	\end{equation}
where the $F_{n+1}^{\mu}$ and $Y_{2(n+1)}^{\mu}$ are irreducible representations of $Gl_{n+1}$ and $S_{2(n+1)}$, respectively, which are associated to Young Diagram $\mu$ with $2(n+1)$ boxes and number of nonzero rows $\leq n+1$. 
We compute the $O_{n+1}$ invariants, 
	\begin{alignat*}{2}
		[\otimes^{2(n+1)}\mathbb{C}^{n+1}]^{O_{n+1}} &\simeq \bigoplus\limits_{\mu} \left[F_{n+1}^{\mu} \otimes Y_{2(n+1)}^{\mu}\right]^{O_{n+1}} \\
		 &\simeq \bigoplus\limits_{\mu} (F_{n+1}^{\mu})^{O_{n+1}} \otimes Y_{2(n+1)}^{\mu}
	\end{alignat*}

By the Cartan-Helgason Theorem, we have that $dim(F_{n+1}^{\mu})^{O_{n+1}}$ is nonzero and equal to one only when the corresponding tableaux $\mu$ has all even parts. Thus, we let $\mu = 2\lambda$:
	\begin{equation}
		\bigoplus\limits_{\substack{{\mu: \hspace{.05cm} \mu = 2\lambda}\\{\mu \vdash 2(n+1)}\\ {\ell(\mu) \leq n+1} }} (F_{n+1}^{\mu})^{O_{n+1}} \otimes Y_{2(n+1)}^{\mu} = \bigoplus\limits_{\substack{{2\lambda: \hspace{.06cm}2\lambda \vdash 2(n+1)} \\ {\ell({2\lambda}) \leq n+1}}} (F_{n+1}^{2\lambda})^{O_{n+1}} \otimes Y_{2(n+1)}^{2\lambda} 
	\end{equation} 	
This space consists of all Young diagrams $2\lambda$, where $2\lambda \vdash 2(n+1)$, and $\ell(2\lambda) \leq n+1$. We know from \cite{Jeb} that this space of invariants corresponds to the set of fixed-point free involutions, thus we have that 
	\begin{equation}\label{schen}
		dim\left(\bigoplus\limits_{\substack{{2\lambda: \hspace{.05cm} 2\lambda \vdash 2(n+1)} \\ {\ell(2\lambda) \leq n+1}}} (F_{n+1}^{2\lambda})^{O_{n+1}} \otimes Y_{2(n+1)}^{2\lambda}\right) = \dfrac{(2(n+1))!}{2^{n+1}(n+1)!}
	\end{equation}
\textbf{Remark:} The Robinson-Schensted correspondence associates a pair of standard Young tableaux, $(P,Q)$ to a permutation. It is shown, \cite{Schen}, that if the permutation is an involution, then $P = Q$. Furthermore, due to a result of Schutzenberger \cite{Schutz}, we have that the fixed-point free involutions correspond to Young diagrams with all even rows. Thus, the space described above in Equation \ref{schen} consists of all standard Young tableaux of shape $2\lambda$, which again have all even rows.

Recall the kernel contained in this space, 
	\begin{equation*}
		\mathcal{R}\mathcal{E}\mathcal{L} \subset \bigoplus\limits_{\substack{{2\lambda: \hspace{.06cm}2\lambda \vdash 2(n+1)} \\ {\ell(2\lambda) \leq n+1}}} (F_{n+1}^{2\lambda})^{O_{n+1}} \otimes Y_{2(n+1)}^{2\lambda}
	\end{equation*}
is the space of relations between the invariants. Once again, we know that relations do not exist in the stable range. Therefore, relations occur when we violate the inequality $\ell(2\lambda) \leq n$. 

As such, $\mathcal{R}\mathcal{E}\mathcal{L}$ is the space of irreducible representations corresponding to the Young diagrams in our space where $\ell(2\lambda) >n$,
	\begin{equation*}
		\mathcal{R}\mathcal{E}\mathcal{L} = \bigoplus\limits_{\substack{ {2\lambda: \hspace{.06cm}2\lambda \vdash 2(n+1)} \\{\ell(2\lambda) > n}     }} (F_{n+1}^{2\lambda})^{O_{n+1}} \otimes Y_{2(n+1)}^{2\lambda} \subset \bigoplus\limits_{\substack{{2\lambda: \hspace{.06cm}2\lambda \vdash 2(n+1)} \\ {\ell(2\lambda) \leq n+1}}} (F_{n+1}^{2\lambda})^{O_{n+1}} \otimes Y_{2(n+1)}^{2\lambda}
	\end{equation*}
since $\bigoplus\limits_{\substack{ {2\lambda \vdash 2(n+1)} \\{\ell(2\lambda) > n}     }} (F_{n+1}^{2\lambda})^{O_{n+1}} \otimes Y_{2(n+1)}^{2\lambda}$ is the kernel of the map:$$
\zeta: [\otimes^{2(n+1)}\mathbb{C}^{n+1}]^{O_{n+1}\times \Delta S_{n+1}} \twoheadrightarrow [\otimes^{2(n+1)}\mathbb{C}^n]^{O_{n} \times \Delta S_{n+1}}.$$
Now, since we have the restriction that the length of $\lambda$ must be greater than $n$, we have only one option for the Young diagram, that is, $\ell(2\lambda) =n+1$,  

\ytableausetup{smalltableaux}
\begin{center}
	\begin{tabular}{r@{}l}

		\raisebox{-3.3ex}{$n+1\left\{\vphantom{\begin{array}{c}~\\[8ex] ~
				\end{array}}\right.$} &
		\begin{ytableau}[]
			&\\ &\\ &\\ &\\ \end{ytableau} \\ 
		\vspace{-.5cm}
		
	\end{tabular}
	$=: 2\lambda$ \\
	\vspace{-1.3cm}
	\hspace{.4cm}\vdots

	\vspace{0cm} \hspace{.27cm}
	\begin{ytableau}[]
		&\\ \end{ytableau}	
\end{center}

Thus we define Young diagram $2\lambda$ as the Young diagram pictured above, with two columns and $n+1$ rows. Here, we have that the dimension of the irreducible representation $Y_{2(n+1)}^{2\lambda}$ is equal to the number of standard Young tableaux of the column shape $(n+1) \times 2$, thus 
	\begin{equation*}\label{catalan}
		dim\left(\bigoplus\limits_{\substack{ { 2\lambda: \hspace{.05cm}2\lambda \vdash 2(n+1)} \\{\ell(2\lambda) > n}     }} (F_{n+1}^{2\lambda})^{O_{n+1}} \otimes Y_{2(n+1)}^{2\lambda}\right) = C_{n+1}
	\end{equation*} 
where $C_{n+1}$ is a Catalan number. 

The irreducible representation $Y_{2(n+1)}^{2\lambda}$ corresponds to partition $[n+1,n+1]$, and we want to restrict this partition to $S_{n+1} \times S_{n+1}$. As dicussed in Section \ref{littlewood}, we can induce representations in the following way: 
	\begin{equation}\label{induce}
		Ind_{S_{n+1} \times S_{n+1}}^{S_{2(n+1)}} Y^{\alpha}_{n+1} \otimes Y^{\beta}_{n+1} = \bigoplus\limits_{2\lambda: \hspace{.06cm} 2\lambda \vdash 2(n+1)} c_{\alpha\beta}^{2\lambda}Y^{2\lambda}_{2(n+1)}
	\end{equation}
Where the coefficients $c_{\alpha\beta}^{2\lambda}$ are the Littlewood-Richardson numbers discussed in Section 3.3. These coefficients count the number of skew semi-standard Young tableaux of shape $2\lambda/\alpha$ with weight $\beta$.

By Frobenius reciprocity for finite groups, we can restate Equation \ref{induce} in terms of restricting the representation,  
	\begin{equation*}
		Res_{S_{n+1}\times S_{n+1}}^{S_{2(n+1)}} Y^{2\lambda}_{2(n+1)} = \bigoplus\limits_{\substack{{\alpha \vdash n+1}\\{\beta \vdash n+1} }} c^{2\lambda}_{\alpha\beta}Y_{n+1}^{\alpha} \otimes Y_{n+1}^{\beta}
	\end{equation*} 
Thus we have the decomposition: 
	\begin{equation*}
		\bigoplus\limits_{\substack{ {2\lambda: \hspace{.06cm} 2\lambda \vdash 2(n+1)} \\{\ell(2\lambda) = n+1}     }} (F_{n+1}^{2\lambda})^{O_{n+1}} \otimes Y_{2(n+1)}^{2\lambda} = \bigoplus\limits_{\substack{{2\lambda: \hspace{.05cm}2\lambda \vdash 2(n+1)} \\{\alpha \vdash n+1}\\{\beta \vdash n+1} }} (F_{n+1}^{2\lambda})^{O_{n+1}} \otimes c^{2\lambda}_{\alpha\beta}(Y_{n+1}^{\alpha} \otimes Y_{n+1}^{\beta})^{\Delta S_{n+1}}
	\end{equation*}	
The Littlewood-Richardson rule tells us $c_{\alpha\beta}^{2\lambda} \neq 0$ when the Young diagrams of $\alpha$ and $\beta$ fit inside the Young diagram of $2\lambda$. Furthermore,  since $|2\lambda| = 2(n+1)$ and $|\alpha| = n+1$ and $|\beta| = n+1$ we must have that $\alpha = \beta$.

Alternately, we determine that $\alpha = \beta$ by considering homomorphisms between the two irreducible representations. We have that as representations, the $S_{n+1}$ are self dual, and thus we can view the tensor of irreducible representations as an endomorphism group: 
	\begin{equation*}
		Y_{n+1}^{\alpha} \otimes Y_{n+1}^{\beta} \cong End(Y_{n+1}^{\alpha}, Y_{n+1}^{\beta})
	\end{equation*}
By Schur's Lemma we see that there are no nonzero homomorphisms between distinct irreducible representations, and thus we must have that $\alpha = \beta$:
	\begin{equation*}
		\bigoplus\limits_{\substack{{2\lambda: \hspace{.05cm}2\lambda \vdash 2(n+1)} \\{\alpha \vdash n+1}\\{\beta \vdash n+1} }} c^{2\lambda}_{\alpha\beta}Y_{n+1}^{\alpha} \otimes Y_{n+1}^{\beta} = \bigoplus\limits_{\substack{{2\lambda: \hspace{.05cm}2\lambda \vdash 2(n+1)}\\{\alpha \vdash n+1} }} c^{2\lambda}_{\alpha\alpha}Y_{n+1}^{\alpha} \otimes Y_{n+1}^{\alpha}
	\end{equation*} 
Now, the Littlewood-Richardson coefficients $c_{\alpha\alpha}^{2\lambda}$ correspond to the number of semi-standard fillings of tableaux of skew-shape $2\lambda/ \alpha$ of weight $\alpha$. Thus, since $2\lambda$ has two columns of length $n+1$, it is clear that we must have that each coefficient is equal to $1$. 

Now, we recall that we are working in the space of relations, $\mathcal{R}\mathcal{E}\mathcal{L}$, and thus we have that 
	\begin{equation*}
		\mathcal{R}\mathcal{E}\mathcal{L}_{n+1} = \bigoplus\limits_{\substack{ {\alpha: \hspace{.06cm} \alpha \vdash n+1} \\ {c_{\alpha \alpha}^{2\lambda}=1} \\{2\lambda \vdash 2(n+1)}\\ {\ell(2\lambda)=n+1}}} (Y_{n+1}^{\alpha} \otimes Y_{n+1}^{\alpha})^{\Delta S_{n+1}}	
	\end{equation*}

So, to each $\alpha$ there exists a polynomial relation which corresponds to the irreducible representation $Y_{n+1}^{\alpha}$. 

The number of $\alpha$ that satisfy this is precisely the dimension of the space of relations,  
	\begin{equation}
		dim(\mathcal{R}\mathcal{E}\mathcal{L}_{n+1}) = \begin{cases}
		\dfrac{n}{2} +1 &n \textnormal{  even} \\
		\dfrac{n+3}{2} &n \textnormal{  odd}
		
		\end{cases}
	\end{equation}
We proceed by induction on $n$, the degree of the invariants. 

Case I: Let $n$ be a positive, even integer; $n = 2m$ for some $m \in \mathbb{N}$. Base case: Let $n=2$. Then we are considering the dimension of the space of relations between the elements of $\mathcal{P}^3(M_2)^{O_2}$. Thus we look at all partitions, $\alpha \vdash 3$. We are concerned with the specific $\alpha$ that give relations in our space; these are the ones in which the Young diagrams corresponding to $\alpha$ fit exactly inside the column shaped Young diagram $(2,2,2)$:
	\begin{figure}[H]
		\begin{center}\resizebox{2cm}{!}{
			\ytableausetup{centertableaux}
			\ydiagram{2,2,2} }
			
			\vspace{.2cm} 
			\footnotesize{Young diagram $(2,2,2)$ }
		\end{center}
	\end{figure}
Below we show all of the possible partitions of $\alpha$ and their corresponding Young diagrams: 
	
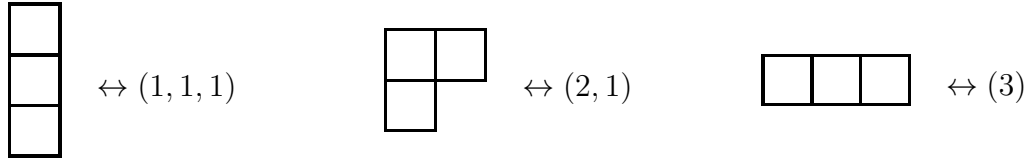
\begin{figure}[H]
		\begin{center}
		\begin{subfigure}{5cm}
			\begin{center}\resizebox{1.25cm}{!}{
			\ytableausetup{centertableaux}
			\ydiagram{1,1,1}} \hspace{.2cm} $\leftrightarrow (1,1,1)$
			\end{center}
		\end{subfigure}
		\begin{subfigure}{5cm}
			\begin{center}\resizebox{1.9cm}{!}{
			\ytableausetup{centertableaux}
			\ydiagram{2,1}} \hspace{.2cm} $\leftrightarrow (2,1)$
			\end{center}
		\end{subfigure}
		\begin{subfigure}{5cm}
			\begin{center}\resizebox{2.5cm}{!}{
			\ytableausetup{centertableaux}
			\ydiagram{3}} \hspace{.2cm} $\leftrightarrow (3)$
			\end{center}
		\end{subfigure}
	\end{center}
	\end{figure} 
Clearly, only partitions $(1,1,1)$ and $(2,1)$ will fit appropriately inside the Young diagram of column shape $(2,2,2)$. Thus, when $n = 2$, we have the dimension of the space of relations is equal to $2 = \dfrac{2}{2} +1$. 

Induction step: We assume the proposition holds for even integer $n = k$, that is, for positive even integer $k$, the dimension of the space of relations between elements of $\mathcal{P}^{k+1}(M_k)^{O_k}$ is equal to $\dfrac{k}{2} +1$. We show this holds for $k+2$. 

Thus, consider the following two partitions $\alpha \vdash k+2+1 = k+3$ and their corresponding Young diagrams: 
	
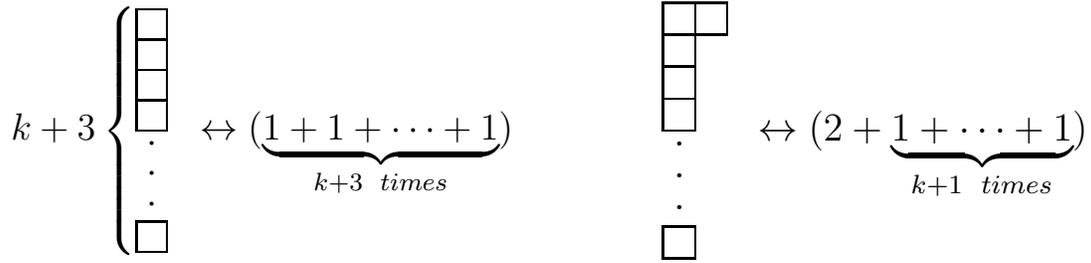
\begin{figure}[H]
	\begin{center}
		\begin{subfigure}{8cm}
			\begin{center}\resizebox{7cm}{!}{
				\ytableausetup{centertableaux}
					 $k+3\left\{\begin{ytableau}[]
					\\ \\ \\ \\ \none[\cdot] \\ \none[\cdot]  \\ \none[\cdot] \\ \\ \end{ytableau}\right. \hspace{.2cm} \leftrightarrow (\underbrace{1 +1 + \dots + 1}_{k+3 \hspace{.2cm} times})$}
			\end{center}
		\end{subfigure}
		\begin{subfigure}{8cm}
			\begin{center}\resizebox{6cm}{!}{
				\ytableausetup{centertableaux}
					 $\begin{ytableau}[]
					&\\ \\ \\ \\ \none[\cdot] \\ \none[\cdot]  \\ \none[\cdot] \\ \\ \end{ytableau} \hspace{.2cm} \leftrightarrow (	2+ \underbrace{1 + \dots + 1}_{k+1 \hspace{.2cm} times})$}
			\end{center}
		\end{subfigure}
	\end{center}
\end{figure} 

Then, it is clear that the Young diagram corresponding to partition $(1 + \dots +1)$ fits inside the diagram $(2, \dots, 2)$ of $2(k+1)$ boxes. Additionally, we know from our assumption that there are $\dfrac{k}{2}+1$ partitions of $k+1$ that satisfy our condition on the corresponding Young diagrams. Thus the total number of partitions that work is 
	\begin{equation*}
		\dfrac{k}{2} +1 +1  = \dfrac{k+2}{2} +1. 
	\end{equation*} 
Case II: Let $n$ be an odd, positive integer; $n = 2m +1$ for some $m \in \mathbb{N}$. 
Base case: Let $n =1$. Then we are considering the dimension of the space of relations between the elements of $\mathcal{P}^2(M_1)^{O_1}$. Thus we look at all partitions, $\alpha \vdash 2$. We are again concerned with only the $\alpha$ which correspond to relations in our space, that is, the $\alpha$ in which their corresponding Young diagrams fit exactly into the column shaped Young diagram $(2,2)$: 
	\begin{figure}[H]
	\begin{center}\resizebox{2cm}{!}{
		\ytableausetup{centertableaux}
		\ydiagram{2,2} }
		
		\vspace{.2cm} 
		\footnotesize{Young diagram $(2,2)$ }
	\end{center}
	\end{figure}
Below we show all the possible partitions of $\alpha$ and the corresponding Young diagrams:
	
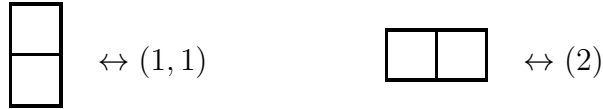
\begin{figure}[H]
	\begin{center}
	\begin{subfigure}{5cm}
		\begin{center}\resizebox{1.25cm}{!}{
		\ytableausetup{centertableaux}
		\ydiagram{1,1}} \hspace{.2cm} $\leftrightarrow (1,1)$
	\end{center}
	\end{subfigure}
	\begin{subfigure}{5cm}
		\begin{center}\resizebox{1.9cm}{!}{
		\ytableausetup{centertableaux}
		\ydiagram{2}} \hspace{.2cm} $\leftrightarrow (2)$
	\end{center}
	\end{subfigure}
	\end{center}
	\end{figure}
Both partitions $(1,1)$ and $(2)$ fit appropriately inside the column shape $(2,2)$. Thus, when $n=1$, we have the dimension of the space of relations is equal to $2 = \dfrac{1+3}{2}$.

Induction step: We suppose the proposition holds for odd integer $n =k$. Therefore we assume for odd integer $k$, the dimension of the space of relations between elements of $\mathcal{P}^{k+1}(M_k)^{O_k}$ is equal to $\dfrac{k+3}{2}$. We show this proposition holds for $k+2$. 

Similar to Case I, we consider the following two partitions of $k+2+1 = k+3$: 
	\begin{equation*}
		\underbrace{1 +1 + \dots + 1}_{k+3 \hspace{.2cm} times}
	\end{equation*}
	\begin{equation*}
		2+ \underbrace{1 + \dots + 1}_{k+1 \hspace{.2cm} times}
	\end{equation*}
Then, we know from our assumption that there are $\dfrac{k+3}{2}$ partitions of the $k+1$ that satisfy our conditions. Thus the total number of partitions that work is
	\begin{equation*}
		\dfrac{k+3}{2} +1 = \dfrac{k+3+2}{2} .
	\end{equation*}
which thus concludes the proof. 

\chapter{Finding a basis of $\mathcal{R}\mathcal{E}\mathcal{L}$ }\label{finding a basis} 
In the previous subsection we state and prove a theorem regarding the dimension of the space $\mathcal{R}\mathcal{E}\mathcal{L}_{n+1}$. Here, we present a method to determine a basis for this space of relations by relying on a construction of elements from the group algebra $\mathbb{C}[S_{2(n+1)}]$, called Young symmetrizers. The proceeding section presents a short explanation of Young symmetrizers, and states several concepts that we will use to obtain a basis of the space of relations. 
\vspace{.1cm}
\newline
\textbf{4.1 A discussion of Young symmetrizers.}\label{YS discussion} We create a standard Young tableau by filling the Young diagram with the numbers $1, \dots, n$ such that the rows and columns strictly increase. Then, we can define the following:
	\begin{alignat*}{2}
		P_{\lambda} &= \{\sigma \in S_{n} | \sigma \textnormal{ preserves each row of } \lambda   \} \\
		Q_{\lambda} &= \{\sigma \in S_{n} | \sigma \textnormal{ preserves each column of } \lambda    \}
\end{alignat*}
These subgroups of $S_n$ define elements $a_{\lambda}$ and $b_{\lambda}$ in the group algebra, $\mathbb{C}[S_{n}]$: 
	\begin{alignat*}{2}
		a_{\lambda} &:= \sum\limits_{\sigma \in P_{\lambda}} e_{\sigma} \\
		b_{\lambda} &:= \sum\limits_{\sigma \in Q_{\lambda}} \textnormal{sgn}(\sigma)e_{\sigma}
	\end{alignat*}	
where $e_{\sigma}$ denotes the unit vector corresponding to $\sigma$. By construction, the elements $a_{\lambda}$ and $b_{\lambda}$ are idempotents in the group algebra. They do not commute, however, their product is also idempotent, and is defined as the \textit{Young Symmetrizer} \cite{Weyl}, \cite{Young}.

\begin{defn}
	The Young symmetrizer corresponding to Young diagram $\lambda$ is defined as \newline $y_{\lambda} := a_{\lambda}\cdot b_{\lambda}$. 
\end{defn} 

We note that any element, $d$, of $\mathbb{C}[S_n]$ gives an invariant subspace, $\mathbb{C}[S_n]d$, of $\mathbb{C}[S_n]$. However, the image of a Young symmetrizer (by right multiplication on $\mathbb{C}[S_n]$) is an invariant subspace which is irreducible under the action of $\mathbb{C}[S_n]$, and unique for each partition $\lambda$ \cite{Fulton}. 

The following theorem tells us that the subspaces $\mathbb{C}[S_n]y_\lambda$ given by the Young symmetrizers are, in fact, irreducible representations of $S_n$, and every irreducible representation of $S_n$ is of this form. 
 
	\begin{thm}\label{young sym thm} Given $S_n$, let $\lambda$ be a partition of $n$. Define $Y^{\lambda}$ as the subspace of $\mathbb{C}[S_n]$ spanned by the Young symmetrizer $y_{\lambda}$. Then:
		\begin{itemize}
			\item $Y^{\lambda}$ is an irreducible representation of $S_n$
			\item If $\lambda, \mu$ are distinct partitions of $n$, then $Y^{\lambda} \ncong Y^{\mu}$
			\item The $Y^{\lambda}$ account for all irreducible representations of $S_n$. 
		\end{itemize}
	\end{thm}
Proof of this theorem can be found in Fulton and Harris \cite{Fulton}. Because we can construct the irreducible representations of $S_n$ in this way, we revisit our decomposition from Section $3$ and construct a basis of $\mathcal{R}\mathcal{E}\mathcal{L}_{n+1}$. 
\vspace{.3cm}
\newline
\textbf{4.2 Constructing a basis using Young symmetrizers.}\label{Constructing a Basis} We show in Section $3$ that the space of relations between invariants under the conjugation action of $O_n$ on $\mathcal{P}(M_n)$, $\mathcal{R}\mathcal{E}\mathcal{L}_{n+1}$ is,
\begin{equation*}
\mathcal{R}\mathcal{E}\mathcal{L}_{n+1} = \bigoplus\limits_{\substack{ {\alpha: \hspace{.06cm} \alpha \vdash n+1} \\ {c_{\alpha \alpha}^{2\lambda}=1} \\{2\lambda: \hspace{.05cm}2\lambda \vdash 2(n+1)}\\ {\ell(2\lambda)=n+1}}} (Y_{n+1}^{\alpha} \otimes Y_{n+1}^{\alpha})^{\Delta S_{n+1}}	
\end{equation*}
where, as vector spaces,  
\begin{equation*}
\bigoplus\limits_{\substack{ {\alpha: \hspace{.06cm} \alpha \vdash n+1} \\ {c_{\alpha \alpha}^{2\lambda}=1} \\{2\lambda \vdash 2(n+1)}\\ {\ell(2\lambda)=n+1}}} (Y_{n+1}^{\alpha} \otimes Y_{n+1}^{\alpha})^{\Delta S_{n+1}} \simeq \bigoplus\limits_{\substack{ {2\lambda:\hspace{.06cm} 2\lambda\vdash 2(n+1)} \\ {\ell(2\lambda)=n+1}    }} (F_{n+1}^{2\lambda})^{O_{n+1}} \otimes Y_{2(n+1)}^{2\lambda}	.
\end{equation*}
For each partition $\alpha$ of $n+1$, there exists a polynomial relation, and this relation corresponds to the $Y_{n+1}^{\alpha}$.

Furthermore, recall that the dimension of the irreducible representation $Y_{2(n+1)}^{2\lambda}$ is equal to a Catalan number, the dimension of the vector space. A natural basis is in one-to-one correspondence with the standard Young tableau with two columns of length $n+1$. These tableaux are all the diagrams 
\ytableausetup{smalltableaux}
\begin{center}
	\begin{tabular}{r@{}l}
		
		\raisebox{-3.5ex}{$n+1\left\{\vphantom{\begin{array}{c}~\\[8ex] ~
				\end{array}}\right.$} &
		\begin{ytableau}[]
			&\\ &\\ &\\ &\\ \end{ytableau} \\ 
		\vspace{-.5cm}
		
	\end{tabular}
	$=: 2\lambda$ \\
	\vspace{-1.2cm}
	\hspace{.4cm}\vdots

	\vspace{0cm} \hspace{.27cm}
	\begin{ytableau}[]
		&\\ \end{ytableau}	
\end{center}
with a standard filling of the numbers $[1, 2, \dots, 2n]$ such that each row and column strictly increases. Let 
\begin{equation*}
T_{2\lambda} := \{ T | \hspace{.1cm} \textnormal{$T$ is a standard tableau of shape $2\lambda$.}  \}.
\end{equation*}
Then we know that each $T$ in the set $T_{2\lambda}$ corresponds to a relation between the elements of $\mathcal{P}^{n+1}(M_n)^{O_n}$. Additionally, since each element of $T_{2\lambda}$ is a standard tableau, it corresponds to an element of $\mathbb{C}[S_{2(n+1)}]$ called a Young symmetrizer. We denote this element $y_{T}$.

Recall in Theorem \ref{my result} we determine the dimension of the space of relations between degree $n+1$ invariants under the conjugation action of $O_n$ on polynomials on $n \times n$ matrices. By averaging over $\Delta S_{n+1}$, we have the projection, 
	\begin{equation*}
		[\otimes^{n+1}(M_n)]^{O_n} \rightarrow \mathcal{P}^{n+1}(M_n)^{O_n}
	\end{equation*}
between the space of orthogonally invariant tensors and the polynomial invariants under the $O_n$ action. We discussed the following maps in Section $2$, 
	\begin{equation*}
		S_{2(n+1)} \rightarrow S_{2(n+1)}/H_{n+1} \rightarrow \Delta S_{n+1}\backslash S_{2(n+1)}/H_{n+1}	.
	\end{equation*}
Recall, the cosets $S_{2(n+1)}/H_{n+1}$ are in one-to-one correspondence with the set $I_{2(n+1)}$, the fixed-point free involutions on $S_{2(n+1)}$. Furthermore, the double cosets in the above map are in bijective correspondence to a basis of the degree $n+1$ polynomial invariants under the orthogonal group action on $\mathcal{P}(M_n)$.

Thus we have the following invariant subspaces of the full group algebra: 
	\begin{equation}\label{group algebra map}
		\mathbb{C}[S_{2(n+1)}]^{\Delta S_{n+1}\times H_{n+1}} \hookrightarrow \mathbb{C}[S_{2(n+1)}]^{H_{n+1}} \hookrightarrow \mathbb{C}[S_{2(n+1)}]
	\end{equation}
where elements of $\mathbb{C}[S_{2(n+1)}]^{\Delta S_{n+1}\times H_{n+1}}$ are linear combinations of the permutations that correspond to the polynomial invariants, $\mathcal{P}^{n+1}(M_n)^{O_n}$. 

Then, we consider the following maps: 
\begin{equation*}\label{epsilon}
	\mathbb{C}[S_{2(n+1)}]^{\Delta S_{n+1} \times H_{n+1}} \xrightarrow{\mathcal{E}} [\otimes^{n+1}M_n]^{O_n \times \Delta S_{n+1}} \rightarrow \mathcal{P}^{n+1}[M_n]^{O_n} 
\end{equation*}
where $\mathcal{E}$ takes elements of the group algebra invariant under the left $\Delta S_{n+1}$ action and the right $H_{n+1}$ action into endomorphisms on tensors. 

So, for arbitrary degree of invariants $d$, we have the projection:
	\begin{equation*}
		\mathbb{C}[S_{2d}]^{\Delta S_{d} \times H_{d}} \twoheadrightarrow \mathcal{P}^d(M_n)^{O_n}
	\end{equation*}
where if $d=n+1$, there exists a nonzero kernel which corresponds to the relations between degree $n+1$ invariants. 

Thus we take a Young symmetrizer and conjugate $\tau = (12)(34)\cdots(2(n+1)-1\hspace{.1cm}2(n+1))$ by each of its terms in order to write it as a linear combination of fixed point free involutions. We can then determine which double coset, $\Delta S_{n+1}\backslash S_{2(n+1)}/H_{n+1}$, each term is in. Thus, we rewrite the Young symmetrizer using coset representatives for each of its terms. We denote this by $\widetilde{y_{T}}$, so that
	\begin{equation*}
		\widetilde{y_{T}} \in \mathbb{C}[S_{2(n+1)}]^{\Delta S_{n+1}\times H_{n+1}}
	\end{equation*}

Now, the $\widetilde{y_{T}}$ form a spanning set for the relations between the degree $n+1$ invariants. In order to find a basis of relations, we look for a subspace of this vector space with dimension dictated by Theorem \ref{my result}. 

We use Python and Sage to write the code for finding relations via this method of Young symmetrizers. Our code runs on the Sagemath cloud with 11 GB disk space, 5 GB of RAM and 1 core. However, the size of the Young symmetrizers grows exponentially and the calculations quickly become too RAM intensive as we increase the degree of the invariants. For example, in $[\otimes^5M_4]^{O_4}$ there are $42$ Young symmetrizers, each with $460,800$ components; the calculation to find the relations using the method discussed here takes a little over $43$ hours to compute on this server. 

In the next chapter, we discuss a much faster method for finding the relations between the invariants under the $O_n$ action on $\mathcal{P}(M_n)$. 

\chapter{A Monte Carlo Method}\label{Monte Carlo}
In this chapter we introduce a new method for determining a basis for the space of relations, $\mathcal{R}\mathcal{E}\mathcal{L}_{n+1}$. We want to avoid the lengthy calculations involved in computing relations using Young symmetrizers, discussed in Section \ref{finding a basis}. We know that the kernel of the map 
	\begin{equation*}
		\zeta: [\otimes^{2(n+1)}\mathbb{C}^{n+1}]^{O_{n+1}\times \Delta S_{n+1}} \twoheadrightarrow [\otimes^{2(n+1)}\mathbb{C}^n]^{O_{n} \times \Delta S_{n+1}}
	\end{equation*}
is exactly our space $\mathcal{R}\mathcal{E}\mathcal{L}_{n+1}$ and consists of relations between the degree $n+1$ invariants under the complex orthogonal group action on $\mathcal{P}(M_n)$. 

We showed in Chapter $3$ that $\mathcal{P}(M_n)^{O_n}$, the algebra of invariant polynomials under the conjugation action of $O_n$, is generated by elements of the form $Tr(x^{a_1}(x^T)^{a_2}x^{a_3}(x^T)^{a_4}\cdots x^{a_M})$
for matrix $x\in M_n$. Products of the above polynomials form a spanning set of the invariants, and they are not linearly independent when the degree of the monomials is greater than $n$. 

We want to find a basis of $\mathcal{R} \mathcal{E}\mathcal{L}$ using elements of the invariant algebra, $\mathcal{P}^{n+1}(M_n)^{O_n}$, and avoid computations with Young symmetrizers. Determining the relations between these polynomials requires solving the nonlinear equations: 
	\begin{equation}
		\sum\limits_{i=1}^{k}y_iTr(x^{a_{i_1}}(x^T)^{a_{i_2}}\cdots) = 0
	\end{equation} 
where $x \in M_n$, $k = dim(\mathcal{P}^{n+1}(M_n)^{O_n})$, (recall that this is the number of double cosets $S_{n+1}\backslash S_{2(n+1)}/H_{n+1}$), and the $y_i$ are constant coefficients in $\mathbb{C}$. 

Solving these nonlinear equations can be computationally challenging. However, we introduce a method for finding relations that allows us to instead solve a $k \times k$ linear system of equations via a Monte Carlo algorithm. 

Each of the invariants $Tr(x^{a_{i_1}}(x^T)^{a_{i_2}}\cdots)$ are constructed using a matrix $x \in M_n(\mathbb{C})$. By definition, 
	\begin{equation*}
		Tr(x^{a_{i_1}}(x^T)^{a_{i_2}}\cdots) \in \mathbb{C}
	\end{equation*}
Therefore, if we randomly generate a matrix $x \in M_n$ we can compute a numerical value for each of the $k$ invariants in $\mathcal{P}^{n+1}(M_n)^{O_n}$. Each equation then becomes linear in the variables $y_i$: 
	\begin{equation*}
		\sum\limits_{i=1}^{k}y_i\underbrace{Tr(x^{a_{i_1}}(x^T)^{a_{i_2}}\cdots)}_{\in \mathbb{C}} = 0
	\end{equation*}
Thus, if we generate $k$ random matrices in $M_n$, and compute the values of each invariant in $\mathcal{P}^{n+1}(M_n)^{O_n}$, this gives $k$ different linear equations. 

Let $x_1, x_2, \dots x_k$ denote the $k$ randomly generated matrices in $M_n(\mathbb{C})$. We can then solve the $k \times k$ linear system for the $y_i$,
	\begin{alignat*}{3}
		y_1Tr(x_1^{a_{1_1}}(x_1^T)^{a_{1_2}}\cdots) + y_2Tr(x_1^{a_{2_1}}(x_1^T)^{a_{2_2}}\cdots) + \dots + y_kTr(x_1^{a_{k_1}}(x_1^T)^{a_{k_2}}\cdots) &=0 \tag{1}\\
		y_1Tr(x_2^{a_{1_1}}(x_2^T)^{a_{1_2}}\cdots) + y_2Tr(x_2^{a_{2_1}}(x_2^T)^{a_{2_2}}\cdots) + \dots + y_kTr(x_2^{a_{k_1}}(x_2^T)^{a_{k_2}}\cdots) &=0 \tag{2}\\
		\vdots \\
		y_1Tr(x_k^{a_{1_1}}(x_k^T)^{a_{1_2}}\cdots) + y_2Tr(x_k^{a_{2_1}}(x_k^T)^{a_{2_2}}\cdots) + \dots + y_kTr(x_k^{a_{k_1}}(x_k^T)^{a_{k_2}}\cdots) &=0 \tag{k}
	\end{alignat*}
and the solution is exactly the relations in the space $\mathcal{R}\mathcal{E}\mathcal{L}_{n+1}$. 

\textbf{Example:} We revisit the space $\mathcal{P}^3(M_2)^{O_2}$ of degree $3$ invariants under the conjugation action of $O_2$ on $M_2$. We know that the dimension of the invariant space is $5$, recall we have the following invariants: 
	\begin{equation*}
	Tr(x^3), \hspace{.2cm} Tr(x^2x^T),\hspace{.2cm} Tr(x^2)Tr(x),\hspace{.2cm} Tr(xx^T)Tr(x),\hspace{.2cm} Tr(x)^3	
	\end{equation*}
We construct a basis for this space using the Monte Carlo method described above. The relations can be described as solutions to the following the nonlinear equations in the invariants: 
	\begin{equation*}
		y_0Tr(x^3) + y_1Tr(x^2x^T) + y_2Tr(x^2)Tr(x) + y_3Tr(xx^T)Tr(x) + y_4Tr(x)^3 = 0
	\end{equation*}
where the $y_i$'s are constant coefficients.

In order to avoid solving a complicated nonlinear system, we first repeatedly generate random complex $2 \times 2$ matrices and compute numerical values for each of the five invariant polynomials. By generating five random matrices, we produce five linear equations using the numerical values for the invariants. That is we generate, 
	\begin{alignat*}{5}
		y_0Tr(x^3) + y_1Tr(x^2x^T) + y_2Tr(x^2)Tr(x) + y_3Tr(xx^T)Tr(x) + y_4Tr(x)^3 &= 0 \tag{1} \\
		y_0Tr(x^3) + y_1Tr(x^2x^T) + y_2Tr(x^2)Tr(x) + y_3Tr(xx^T)Tr(x) + y_4Tr(x)^3 &= 0 \tag{2} \\	
		\vdots \\
		y_0Tr(x^3) + y_1Tr(x^2x^T) + y_2Tr(x^2)Tr(x) + y_3Tr(xx^T)Tr(x) + y_4Tr(x)^3 &= 0 \tag{5}	
	\end{alignat*}
where each equation uses a different randomly generated $2 \times 2$ matrix to compute a numerical value for each of the invariants, and thus the equations are all linear in the $y_i$ with coefficients in $\mathbb{C}$. 

We then run code to solve this $5 \times 5$ linear system, and the result is two linearly independent relations in the space $\mathcal{R}\mathcal{E}\mathcal{L}$, 
	\begin{alignat*}{2}
		Tr(x)^3 -3Tr(x^2)Tr(x) +2Tr(x^3) &= 0 \\
		2Tr(x^2x^T) - 2Tr(xx^T)Tr(x) + Tr(x)^3 - Tr(x^2)Tr(x) &= 0
	\end{alignat*} 

Again, we use Python and Sage to write the code for finding relations via this Monte Carlo simulation method. Our code runs on the Sagemath cloud with 11 GB disk space, 5 GB of RAM and 1 core. However, this method is much faster at finding relations than the method of Young symmetrizers. For example, in the case of $[\otimes^5M_4]^{O_4}$ which took over $43$ hours to compute relations via Young symmetrizers, our new method takes just a bit over $2$ minutes.

\footnotesize{Alison Elaine Becker, Email: alison.elaine90@gmail.com}

    \end{document}